\input amstex
\documentstyle{amsppt}

\input psfig
\magnification=1200
%\pagewidth{5.4in}
%\pageheight{7.5in}
%\parskip 5pt
%\expandafter\redefine\csname logo\string@\endcsname{}
\NoBlackBoxes
\NoRunningHeads

\topmatter
\title Finiteness of simple homotopy type up to $s$-cobordism of
aspherical $4$-manifolds\\
\endtitle
\thanks {\it 2000 Mathematics Subject Classification.} {\it Primary.} 
57N13, {\it Secondary.} 57R67, 14R10.
\endthanks
\thanks {\it Keywords and phrases.} Structure set, topological rigidity, 
surgery exact sequence, cobordism, $4$-manifolds, Haken
manifolds, Whitehead group, assembly map.\endthanks
\thanks accepted for publication in the Internat. Math. Res.
Notices.\endthanks
\author Sayed K. Roushon
\endauthor
\address School of Mathematics, Tata Institute of Fundamental 
Research, Homi Bhabha Road, Mumbai 400 005, India. 
\endaddress
\email roushon\@math.tifr.res.in \endemail
%\date October 16, 2000 \enddate
\abstract
In this paper we show that for a large class $\Cal C$ of $4$-manifolds 
each member of $\Cal C$ has only finitely many simple homotopy type 
up to $s$-cobordism. This result generalizes a similar 
result of Hillman for certain complex surfaces. We also present a
correction in the proof of Hillman's result.
\endabstract
%\keywords 
%\endkeywords
%\subjclass 
%\endsubjclass
\endtopmatter
\document
\baselineskip 14pt

\head
0. Introduction
\endhead
The Borel conjecture says that closed aspherical manifolds  
are determined by their fundamental groups, i.e., an isomorphism between
the fundamental groups of two closed aspherical manifolds is induced by a
homeomorphism of the manifolds. In dimension greater than 
$4$ this question is answered in positive for the class of manifolds 
with nonpositively curved Riemannian metric: this is the largest class 
of manifolds for which the answer is known so far. The advantage in higher 
dimension is the availability of the $s$-cobordism theorem and the 
surgery theory. In dimension $3$ the answer is known for a vast class 
of manifolds: namely, Haken manifolds and hyperbolic manifolds. The
answer
will be complete in dimension $3$ provided Thurston's Geometrization
conjecture is true. The dimension $4$ case is also not yet settled. For
example the $s$-cobordism theorem and the exactness of the surgery
sequence is known only for $4$-manifolds with elementary amenable
fundamental groups. 

In this paper we show that for a class of aspherical $4$-manifolds; 
up to $s$-cobordism there are only finitely many $4$-manifolds
simple homotopy equivalent to a given member of this class. Due to the
unavailability 
of the $4$-dimensional $s$-cobordism theorem we cannot quite say that
there are only finitely 
many simple homotopy type up to homeomorphism of a manifold from this
class.

\head
1. Finiteness of simple homotopy type
\endhead
Let $\Cal C$ be the class of compact orientable aspherical $4$-manifolds
so that for each $M\in {\Cal C}$ there is a fiber bundle projection 
$M\to {\Bbb S}^1$ with irreducible fiber $N$ and either $M$ has
nonempty boundary or it is closed and satisfies one of the
following properties:
\roster
\item $H^1(N, {\Bbb Z})\neq 0$
\item there is an embedded incompressible torus $T$ in $N$ so that 
$\pi_1(T)$ has square root closed image in $\pi_1(N)$
\item $N$ is a Seifert fibered space with base surface of genus $\geq 1$
\item $N$ has more than $2$ geometric pieces in its Jaco-Shalen and
Johannson decomposition and the dual graph of this decomposition has a
vertex which disconnects the graph and also the fundamental group of any
edge emanating from this vertex is square root closed in the
fundamental group of the target vertex
\item $N$ supports a  hyperbolic metric.
\endroster

Here recall that a subgroup $H$ of a group $G$ is called {\it
square 
root closed} if for any $x\in G$, $x^2\in H$ implies $x\in H$. And a 
$3$-manifold is called {\it irreducible} if any embedded $2$-sphere 
in it bounds an embedded $3$-disc. An irreducible $3$-manifold with 
nonempty boundary has nonvanishing first Betti number. Also any 
irreducible $3$-manifold with nonvanishing first Betti number 
is Haken. 

The dual graph of the $JSJ$-decomposition has vertices the pieces in the
decomposition and edges are
tori which are common boundary component of two pieces. By fundamental
group of a vertex or an edge we mean the fundamental group of the
associated spaces.

Before we state our main theorem we recall the definition of
homotopy-cobordant structure 
sets: Let $M$ be a compact manifold. Define ${\Cal S}^s_{TOP}(M, \partial
M)=\{(N, f)\ |\ f:N\rightarrow M,$ where $N$ a compact
manifold, $f$ a simple homotopy equivalence, $f|_{\partial N}$ is a 
homeomorphism onto $\partial M$ $\}/\simeq$, where $(N_1, f_1)\simeq (N_2,
f_2)$ if there is
a map $F:W\rightarrow M$ with domain $W$ a $s$-cobordism with 
$\partial W=N_1\cup N_2$ and $F|_{N_i}=f_i$. If the Whitehead group 
of $\pi_1(M)$ vanishes then this is the usual homotopy-topological
structure set of $M$ provided the $s$-cobordism theorem is true 
in dim$M$. The dimension $4$ $s$-cobordism theorem is known to be true
only for $4$-manifolds with elementary amenable fundamental group 
(see \cite{FQ}).    

In this paper we prove the following theorem:

\proclaim{Theorem 1.1} Let $M\in {\Cal C}$. Then the set
${\Cal S}^s_{TOP}(M)$ is finite when $M$ is closed and for 
$n\geq 1$ ${\Cal S}^s_{TOP}(M\times {\Bbb D}^n, \partial (M\times {\Bbb
D}^n))$ has only one element.
\endproclaim

\proclaim{Corollary 1.2} Let $M\in {\Cal C}$ and $N$ be any other
$4$-manifold homotopy equivalent to $M$. Then there are integers 
$r$ and $s$ so that $M\#r({\Bbb S}^2\times {\Bbb S}^2)$ is 
diffeomorphic to $N\#s({\Bbb S}^2\times {\Bbb S}^2)$. Here $\#$ denotes
connected sum. In such a case $M$ and $N$ are called stably
diffeomorphic.\endproclaim

\demo{Proof of Corollary 1.2} Follows from the main theorem in
\cite{D}.\enddemo

\remark{Remark 1.3} As it is not yet known if a $s$-cobordism between two
$4$-manifolds is trivial we cannot quite conclude
that the manifolds in the class $\Cal C$ has finitely many homotopy 
type up to homeomorphism.\endremark

Here we deduce an interesting corollary: 

\proclaim{Corollary 1.4} Let $M$ be a nonsingular complex affine surface
(i.e., a nonsingular complex algebraic surface in the 
complex space ${\Bbb C}^n$) which is a fiber bundle over the
circle with irreducible (in $3$-manifold sense) fiber. Then for $n\geq 1$
$M\times {\Bbb D}^n$ has only one homotopy type (with homotopy which are
homeomorphism outside a compact set) up to homeomorphism.\endproclaim

\demo{Proof} Using a suitable Morse function on $M$ (for example 
consider the polynomial function $||x-x_0||^2$ for a fixed $x_0\in M$) it
is easily deduced (by Morse theory) that $M$ is diffeomorphic to the
interior of a compact 
aspherical $4$-manifold. This follows because the restriction to $M$ 
of any polynomial function has only finitely many critical value. 
(see corollary 2.8 in \cite{M}). The Corollary follows.\qed\enddemo    

\head
2. Proof of the theorem 1.1
\endhead

At first we check that the fundamental group of any
of the $4$-manifolds in the class $\Cal C$ has vanishing Whitehead 
group. 

If $\partial N$ is nonempty and has ${\Bbb S}^2$ as a
boundary component then by irreducibility $N$ is homeomorphic to ${\Bbb
D}^3$ and hence $M$ is homeomorphic to 
${\Bbb D}^3\times {\Bbb S}^1$. In this particular case the
theorem is known. So we assume that if $\partial N\neq \emptyset$ then
genus of any component of $\partial N$ is $\geq 1$. 

Note that for any $M\in {\Cal C}$ the fiber $N$ of the fiber bundle 
$M\rightarrow {\Bbb S}^1$ is a Haken $3$-manifold in the cases
$(1)-(4)$. This implies
$\pi_1(N)\in {\Cal Cl}$ from the notation of \cite{W} and hence 
it has vanishing Whitehead group. Also $\pi_1(N)$ is regular coherent. 

Now we can use the Mayer-Vietoris exact sequence (for $K$-theory) 
from \cite{W} (Sections 17.1.3 and 17.2.3) to deduce that $\pi_1(M)$ has
vanishing Whitehead group. 

A general version of this fact is proved in (lemma V.3, \cite{H1}).

If $N$ is hyperbolic then by the Mostow rigidity theorem the monodromy 
diffeomorphism of the fiber bundle $M\rightarrow {\Bbb S}^1$ is homotopic to
an isometry of finite order and hence $\pi_1(M)$ is isomorphic to the
fundamental group of a $4$-manifold $M'$ which has a ${\Bbb H}^3\times
{\Bbb R}$ structure. Since $M'$ is nonpositively curved
$Wh(\pi_1(M))=Wh(\pi_1(M'))=0$ by \cite{FJ}. 
This conclusion also can be made by noting that in fact the monodromy 
diffeomorphism is (topologically) isotopic (see \cite{G} and \cite{GMT}) 
to a (finite order) isometry and hence the fiber bundle $M$ itself has
${\Bbb H}^3\times {\Bbb R}$ structure.  

In \cite{Ro1} and \cite{Ro2} we proved the following theorem: 

\proclaim{Theorem 2.1} (Theorem 1.2 in \cite{Ro1} and 
Theorem 1.1 and 1.3 in
\cite{Ro2}) Let $N$ be a compact orientable irreducible $3$-manifold 
so that one of the following properties is satisfied:
\roster
\item $H^1(N, {\Bbb Z})\neq 0$
\item there is an embedded incompressible torus $T$ in $N$ so that 
the image of $\pi_1(T)$ is square root closed in $\pi_1(N)$
\item $N$ has more than $2$ geometric pieces in its Jaco-Shalen and
Johannson decomposition and the dual graph of this decomposition has a
vertex which disconnects the graph and also the fundamental group of any
edge emanating from this vertex is square root closed in the
fundamental group of the target vertex.\endroster

Then  for $n\geq 2$ $N\times {\Bbb D}^n$ has only one homotopy type up to 
homeomorphism. 
\endproclaim

Here note that the case when $N$ has nonempty boundary is included in 
case $(1)$. In \cite{Ro2} a large class of  examples of $3$-manifolds is
given satisfying the property $(2)$ and $(3)$ in the above theorem. In
fact it was 
shown there that if we consider the Jaco-Shalen and Johannson ($JSJ$) 
decomposition of the Haken manifold with $T$ as one of the decomposing 
torus then the square root closed condition depends only on the pieces
which abut the torus $T$. Also a large class of examples of
$3$-manifolds are given which has a square root closed incompressible
torus boundary component. 

We recall the Wall-Novikov surgery exact sequence here:

Let ${\Cal S}(X, \partial X)$ denote the topological structure set of $X$
of the group
of homotopy type of $X$ up to homeomorphism. For precise
definition see any reference on surgery theory or in \cite{Ro1}. (Here
note that the differentiable structure set is not a group.) In terms
of this group Theorem 2.1 say that ${\Cal S}(N\times {\Bbb D}^n,
\partial({\Cal S}(N\times {\Bbb D}^n))$ is
trivial for $n\geq 2$. We always assume that $Wh(\pi_1(X))=0$. Then these
groups fit into a long exact sequence of groups:  
$$\cdots\longrightarrow
{\Cal S}_{n-1}(X)\longrightarrow H_n(X, {\Bbb L}_0)\longrightarrow
L_n(\pi_1(X))\longrightarrow {\Cal S}_n(X)\longrightarrow\cdots$$

Here ${\Cal S}_n(X)$ are the total surgery obstruction group of Ranicki
(\cite{R1}) and they are in bijection with ${\Cal S}(X\times {\Bbb D}^n,
\partial(X\times {\Bbb D}^n))$ with a different indexing. 

Note that the fundamental group of any $M\in {\Cal C}$ is of the form 
$\pi_1(M)=\pi_1(N)\rtimes {\Bbb Z}$. From Theorem 2.1, the main theorem in
\cite{S} for the case $(3)$ and by
Farrell-Jones Topological Rigidity theorem for nonpositively curved
Riemannian manifold (in the hyperbolic case) (\cite{FJ}) it follows that
the (assembly) 
map $H_n(N, {\Bbb L}_0)\longrightarrow L_n(\pi_1(N))$ is an isomorphism for
large $n$. 

Now the proof of the theorem follows from the following facts:

{\bf Fact 1}: the Whitehead group of $\pi_1(M)$ vanishes.

{\bf Fact 2}: the Ranicki Mayer-Vietoris exact sequence of surgery 
groups for groups which are semidirect product of a group with the 
infinite cyclic group (see \cite{R2}).

{\bf Fact 3}: the Mayer-Vietoris exact sequence of the generalized 
homology theory for $K(\pi, 1)$ spaces with coefficient in the surgery
spectrum ${\Bbb L}_0$.

{\bf Fact 4}: naturality of the assembly map and an application
of five-lemma together with Theorem 2.1 and Siebenmann's periodicity
theorem (\cite{KS}).

{\bf Fact 5}: the corollary to the theorem V.12 in \cite{H1} which says
that if the assembly map $H_5(M, {\Bbb L}_0)\to L_5(\pi_1(M))$ is an
epimorphism then the set ${\Cal S}^s_{TOP}(M)$ is finite.
\qed

\remark{Remark 2.2} Here we remark that the Theorem 1.1 will be true 
for all compact $4$-manifolds which fiber over the circle if Thurston's  
conjecture is true: i.e., if any aspherical closed $3$-manifold is either
Haken, hyperbolic or Seifert fibered space, and if Theorem 2.1
is true for any Haken $3$-manifold.\endremark

\remark{Remark 2.3} Hillman informed the author that he thinks it
follows from \cite{H3} that if $M'$ is homotopically equivalent to 
$M\in {\Cal C}$ then in fact $M'$ and $M$ are $s$-cobordant.\endremark 

\remark{Remark 2.4} In \cite{H1} Hillman proved Theorem 1.1  for the case
when $M$ also supports a complex structure. Also note that in the 
theorem the case of certain complex surfaces is included in case $(3)$
because in (\cite{H2}) it was proved if a complex surface fibers
over the circle then the fiber is a Seifert fibered space. Hillman proved 
the Theorem 1.1 when the fiber is an arbitrary Seifert fibered space
$N$ assuming that $N$ supports a nonpositively curved Riemannian metric 
(page 81 in \cite{H1}). However, the unit tangent bundle of a closed
oriented surface of genus $\geq 2$ is a Seifert fibered space which
has $\widetilde {{\Bbb {SL}}(2, {\Bbb R})}$ structure but does
not support any nonpositively curved Riemannian metric (see \cite{Ro1}). 
\endremark

\medskip
\noindent
{\bf Acknowledgement:} The author would like to thank J.A. Hillman for
sending his reprints. He also thanks Mehta Research Institute in
Allahabad where this work was done.


\newpage
\Refs

\widestnumber\key{\bf Ro1}
\widestnumber\key{\bf KS}
\widestnumber\key{\bf Ro2}
\widestnumber\key{\bf GMT}

\ref\key{\bf C1}
\by S. E. Cappell
\paper A splitting theorem for manifolds
\jour Invent. Math.
\vol 33 
\yr 1976 
\pages 69--170
\endref

\ref\key{\bf C2}
\bysame  
\paper Mayer-Vietoris sequence in Hermitian K-theory 
\inbook Algebraic K-theory III, Lecture Notes in Math. 343
\publ Springer-Verlag
\publaddr New York
\yr 1973
\pages 478--512
\endref

\ref\key{\bf D}
\by J. F. Davis
\paper The Borel/Novikov conjectures and stable diffeomorphism of
$4$-manifolds
\jour preprint, Indiana University
\endref

\ref\key{\bf FJ}
\by F. T. Farrell and L. E. Jones
\paper Topological rigidity for 
compact nonpositively curved manifolds
\jour Proc. Sympos. Pure Math.
\vol 54 Part 3 
\publ Amer. Math. Soc., Providence, R.I. 
\yr 1993
\pages 229--274
\endref

\ref\key{\bf FQ}
\by M. H. Freedman and F. Quinn 
\book Topology of $4$-manifolds
\publ Princeton University Press
\publaddr Princeton 
\yr 1990
\endref

\ref\key{\bf G}
\by D. Gabai
\paper On the geometric and topological rigidity of hyperbolic
$3$-manifolds 
\jour J. Amer. Math. Soc.
\vol 10 
\yr 1997 
\pages 37--74 
\endref

\ref\key{\bf GMT}
\by D. Gabai, G. R. Meyerhoff, N. Thurston 
\paper Homotopy Hyperbolic 3-Manifolds are Hyperbolic
\jour math.GT/9609207
\yr 1996 
\endref

\ref\key{\bf H1}
\by J.A. Hillman
\book The algebraic characterization of geometric $4$-manifolds
\publ Cambridge University press
\publaddr Press Syndicate of the University of Cambridge, Cambridge, New
York, Melbourne
\yr 1994
\endref

\ref\key{\bf H2}
\bysame 
\paper On $4$-dimensional mapping tori and product geometries
\jour J. London Math. Soc. (2)
\vol 58
\yr 1998
\pages 229--238
\endref

\ref\key{\bf H3}
\bysame
\paper Unpublished manuscript
\endref

\ref\key{\bf KS}
\by R. C. Kirby and L. Siebenmann
\book Foundational Essays on Topological Manifolds, 
Smoothings, and Triangulations, Ann. of Math. Studies
\publ Princeton University Press
\publaddr Princeton
\yr 1977
\endref

\ref\key{\bf M}
\by J. Milnor
\book Singular points of complex hypersurfaces, Ann. of Math.
Studies 
\publ Princeton University Press 
\publaddr Princeton
\yr 1968 
\endref

\ref\key{\bf R1}
\by A. Ranicki
\paper The total surgery obstruction
\inbook Lecture Notes in Math. 763
\publ Springer-Verlag
\publaddr New York 
\yr 1979
\pages 275--316
\endref

\ref\key{\bf R2}
\bysame
\paper Algebraic $L$-theory III, Twisted Laurent extension
\inbook Algebraic $K$-theory, III, Lecture Notes in Math. 343
\publ Springer-Verlag
\publaddr Berlin and New York
\yr 1973
\pages 412--463
\endref

\ref\key{\bf Ro1}
\by S.K. Roushon
\paper $L$-theory of $3$-manifolds with non-vanishing
first Betti number
\jour Internat. Math. Res.
Notices
\vol 2000, no. 3.
\pages 107--113 
\endref

\ref\key{\bf Ro2}
\bysame
\paper Vanishing structure set of Haken $3$-manifolds 
\jour Math. Ann.  
\vol 318
\yr 2000 
\pages 609--620
\endref

\ref\key{\bf S}
\by C. W. Stark
\paper Structure set vanish for certain bundles over 
Seifert manifolds
\jour Trans. of Amer. Math. Soc. 
\vol 285 
\yr 1984
\pages 603--615
\endref

\ref\key{\bf W}
\by F. Waldhausen 
\paper Algebraic K-theory of generalized free 
products, Parts 1 and 2
\jour Ann. of Math. 
\vol 108 
\yr 1978
\pages 135--256
\endref

\endRefs
\enddocument